\documentclass[12pt]{amsart}
\usepackage{amsmath,amssymb, euscript, graphicx, hyperref}

\newtheorem{thm}{Theorem}[section]
\newtheorem{lem}[thm]{Lemma}
\newtheorem{prop}[thm]{Proposition}
\theoremstyle{definition}
\newtheorem{defn}[thm]{Definition}

\newtheorem{rem}[thm]{Remark}
\newtheorem{cor}[thm]{Corollary}

\newcommand{\Diff}{\mathop{\rm Diff}}
\newcommand{\Homeo}{\mathop{\rm Homeo}}

\begin{document}

\author{Azer Akhmedov}
\address{Azer Akhmedov,
Department of Mathematics,
North Dakota State University,
PO Box 6050,
Fargo, ND, 58108-6050}
\email{azer.akhmedov@ndsu.edu}

\title[Finiteness of the topological rank of diffeomorphism groups]{Finiteness of the topological rank of diffeomorphism groups}

\begin{abstract}  For a compact smooth manifold $M$ (with boundary) we prove that the topological rank of the diffeomorphism group $\Diff^k_0(M) \ (\Diff^k_0(M,\partial M))$ is finite for all $k\geq 1$. This extends a result from \cite{AC} where the same claim is proved in the special case of $dim M = k = 1$.
\end{abstract}

\maketitle


\section{Introduction}

The {\em topological rank} of a separable topological group $G$, denoted by $\mathrm{rank}_{top}(G)$, is the minimal integer $n\geq 1$ such that for some $n$-tuple $(g_1, \dots , g_n)$ of elements of $G$ the group $\Gamma=\langle g_1,g_2,...,g_n\rangle$ which they generate is dense in $G$. In this case we say $G$ is \textit{topologically $n$-generated} and the $n$-tuple $(g_1, \dots , g_n)$ {\em topologically generates} $G$. 

 It is not difficult to see that the topological tank of $\mathbb{R}^n$ equals $n+1$ while the topological rank of the $n$-torus $\mathbb{T}^n$ equals 1. It is also immediate to see that a $1$-generated group is necessarily Abelian.

 Since $\mathrm{rank}_{top}(\mathbb{T}^n) = 1$, all compact connected solvable Lie groups have topological rank equal to 1. On the other hand, it has been proved by Kuranishi \cite{K} that all compact semisimple Lie groups are 2-generated. Hofmann and Morris \cite{HM} have extended this result to all separable compact connected groups. Very recently, applying the solution of Hilbert's Fifth Problem, Gelander and Le Maitre \cite{GL} have shown that every separable connected locally compact group is topologically finitely generated.
 
 In the realm of more general Polish groups, interesting finiteness results have been obtained by Macpherson in \cite{Mac}, and by A. Kechris and C. Rosendal in \cite{KR}, for the topological ranks of the following groups: $S_{\infty }$ (the group of bijections of $\mathbb{N}$), the automorphism group of the (so-called) random graph (see \cite{Mac}), $H(2^{\mathbb{N}})$ (the homeomorphism group of the Cantor space), $H(2^{\mathbb{N}}, \sigma )$ (the group of the measure preserving homeomorphisms of $2^{\mathbb{N}}$ with the usual product measure), Aut$(\mathbb{N}^{<\mathbb{N}})$ (the automorphism group of the infinitely splitting rooted tree), and Aut$([0,1],\lambda )$ (Lebesgue measure preserving automorphisms of the closed interval), see \cite{KR} (as noted in \cite{KR}, that the latter group is topologically 2-generated is first shown in \cite{Gr} and \cite{P} by different means). In all these examples, it is shown that the groups in question admit a {\em cyclically dense conjugacy class} \footnote{a group $G$ is said to admit a cyclically dense conjugacy class if for some $f, g\in G$, the set $\{f^ngf^{-n} : n\in \mathbb{Z}\}$ is dense in $G$.} thus the topological rank equals two. In particular, all these groups have a {\em Rokhlin property}, i.e. they possess a dense conjugacy class.   

 In \cite{AC}, the finiteness of the ranks of the diffeomorphisms groups of compact 1-manifolds has been proved. These diffeomorphism groups are quite different in many regards from the groups mentioned in the previous paragraph, in fact, one can almost view them as Polish groups at the other end of the spectrum. The current paper can be viewed as a continuation of \cite{AC}; here, we establish finiteness results for the topological rank of the diffeomorphism group $\Diff^k_0(M)$ (the connected component of the group of all diffeomorphisms of $M$) for an arbitrary compact smooth manifold $M$, as well as for the groups $\Diff^k_0(M, \partial M)$ for an arbitrary compact smooth manifold $M$ with boundary. Here, $\Diff^k_0(M, \partial M)$ denotes the connected component of the group of all diffeomorphisms of $M$ fixing the boundary $\partial M$ pointwise. We are also motivated to understand discrete and dense finitely generated subgroups of these groups, seeking parallels with the very rich and mostly established theory of discrete and dense subgroups of Lie groups. The reader may consult the paper \cite{A} for a large number of remarks and open questions in this program where mainly the group $\Diff_+(I)$ has been considered, however, quite many of the questions there are still meaningful for a general compact smooth manifold and in the higher regularity as well.

  We would like to emphasize that the easier case of a finiteness of the rank of $\Homeo_0(M)$ for compact topological manifolds (with boundary) also follows from the proof, although this is a much weaker result. In fact, it is quite well known that the group $\Homeo_+(I)$ of orientation-preserving homeomorphisms of the interval $I=[0,1]$ is topologically $2$-generated, as the Thompson's group $F$ embeds there densely. The standard faithful representation of $F$  can be modified to obtain $C^{\infty }$ smooth representation \cite{GS},  however, this smoothed out representations is not dense in the $C^1$-metric, quite the opposite, it is in fact $C^1$-discrete.
  
  The homeomorphism groups $\Homeo_0(M)$ resemble compact spaces in many ways (especially for the cases of interval $I$, the disc $D^2$ and the 2-sphere $\mathbb{S}^2$), while already for the the case of the interval $M = I$, and in the regularity $k=1$, the group $\Diff_+(I)$ is homeomorphic to an infinite dimensional Banach space. \footnote {Any two infinite dimensional separable Banach spaces are indeed homeomorphic by a result of M.I.Kadets \cite{Ka}; one can also establish an explicit homeomorphism by the map $f\to \ln f'(t)- \ln f'(0)$ from $\Diff_+^1(I)$ to the Wiener space $C_0[0,1]$ - the space of continuous functions $f:[0,1]\to \mathbb{R}$ with $f(0) = 0$.}

 The following theorem is proved in \cite{AC} in the case of $k=1$.

\begin{thm} \label{thm:main} For each $1\leq k\leq\infty$, there exists a finitely generated dense subgroup of $\Diff_+^k(I)$.
\end{thm}

 Our construction relies on the perfectness results of Tsuboi \cite{T} for diffeomorphism groups of the interval, as well as the following lemma on approximation by diffeomorphisms with iterative $n$-th roots which is independently interesting.

\begin{lem}\label{lem:introroot} Let $f\in\Diff_+^k(I)$ without a fixed point in (0,1) and assume that $f$ is $C^k$-tangent to the identity at $0$ and $1$.  Then for every $r>0$, there exists a $g\in\Diff_+^k(I)$ and a positive integer $N$ such that $g$ is $r$-close to identity and $g^N$ is $r$-close to $f$, in the $C^k$ metric.
\end{lem}

 Lemma \ref{lem:introroot} is an extension of the Lemma 1.2 from \cite{AC} to the higher regularity. The proof of this lemma somewhat follows the one of Lemma 1.2 from \cite{AC} but it turns out to be much more involved. As for the perfectness result, we use a more subtle version of it from \cite{T}. 
 
 By the standard fragmentation argument the following corollary is immediate.
  
\begin{cor}  For each $1\leq k\leq\infty$, there exists a finitely generated dense subgroup of $\Diff_+^k(\mathbb{S}^1)$.
\end{cor}

 The main construction in the proof of Theorem \ref{thm:main} generalizes to the case of an arbitrary (smooth) manifold. We then use an analogue of Tsuboi's perfectness result for the diffeomorphism groups of arbitrary compact (smooth) manifolds, recently re-proven by S.Haller, T.Rybicki and J.Teichmann, \cite{HRT}. Then we find a suitable extension of the Lemma \ref{lem:introroot} in the general manifold case thus obtaining the following more general
 
 \begin{thm}\label{thm:manifold} Let $M$ be a compact smooth manifold (with boundary). Then for all $k\geq 1$ the group $\Diff^k_0(M) \ (\Diff^k_0(M,\partial M))$ admits a finitely generated dense subgroup.
 \end{thm}   
 
 It should be noted that higher dimensional perfectness results for diffeomorphism groups appear already in the works of M.R.Herman \cite{He1, He2} and W.Thurston \cite{Th}. The perfectness of the diffeomorphism group $\Diff _{c}^{\infty }(M)_0$ for a smooth manifold $M$ has been proved by D.Epstein in \cite{Ep} which in turn uses the ideas of J.Mather in \cite{Ma1}, \cite{Ma2} who dealt with the $C^k$ case, $1\leq k < \infty , \ k\neq n+1$. In these works, for an arbitrary diffeomorphism $g\in \Diff _{c}^{\infty }(M)$, a representation $g = [u_1, v_1]\dots [u_N, v_N]$ is guaranteed but seems with less control on $u_i, v_i$, and with a weaker bound on $N$. In \cite{HRT}, besides establishing the {\em uniform perfectness} with a very strong bound, also the so-called {\em smooth perfectness} of the group $\Diff _{c}^{\infty }(M)_0$ has been proved.

 \vspace{1cm}
 
 \section{Finiteness of the rank in higher regularity and in higher dimensions}
 
 In this section, we will prove Theorem \ref{thm:main} and Theorem \ref{thm:manifold} modulo the axillary results Lemma \ref{lem:introroot} and the extension of it in the $n$-dimensional case, Lemma \ref{lem:rootextlem}. The general scheme of the major construction follows the one of \cite{AC}.   

 \medskip
 
  We would like to make an important observation that the proof of Theorem 2.1 in \cite{AC} works in showing that Diff$^k_{+}(I)$ is finitely generated for all $k\geq 1$. Indeed, in the proof of Theorem 2.1, one only needs to make the following modifications: 
  
   \medskip
  
  First, the group $G$ is redefined to be equal to $$\{f\in \mathrm{Diff}_{+}^k(I) \ | \ f \ \mathrm{is} \ C^k\mathrm{-flat}\}$$ where $C^k$-{\em flatness} means that $f'(0) = f'(1) = 1 \ \mathrm{and} \  f^{(j)}(0) = g^{(j)}(0) = 0, 2\leq j\leq k.$ Secondly, the $C^1$-metric everywhere is replaced with $C^k$-metric, and the $C^1$-spaces are replaced with $C^k$-spaces. 
  
  \medskip
  
  Thirdly, the condition (b-i) is replaced with 
  
  (b-i)$' \ \theta '(a_1) = \theta '(b_1)  = 1$ and $\theta ^{(j)}(a_1) = \theta ^{(j)}(b_1) = 0$ for all $2\leq j\leq k$ where $\theta $ is any of the functions $f, g, u$ and $v$. 
  
  Also, the condition (c-i) is replaced with the following condition:
  
  (c-i)$'$ for all $x\in \{a_1, \dots , a_n\}\cup \{b_1, \dots , b_n\}$, $f'(x) = g'(x) = u'(x) = v'(x) = 1$ and $f^{(j)}(x) = g^{(j)}(x) = u^{(j)}(x) = u^{(j)}(x) = 0, 2\leq j\leq k$;
  
  \medskip
  
  and the condition (d-i) is replaced with
  
   (d-i)$'$ $f'(x) = g'(x) = u'(x) = v'(x) = 1$ and $f^{(j)}(x) = g^{(j)}(x) = u^{(j)}(x) = v^{(j)}(x) = 0, 2\leq j\leq k$, for all $x\in \{a_{n+1}, b_{n+1}\}$.
   
   \medskip
   
   Besides, we also assume that the function $h$ is $C^{k}$-flat (instead of being $C^1$-flat) at the endpoints of the intervals $J_n, n\geq 1$.
    
   \medskip
   
   In the main body of the construction, instead of the simple version of the Tsuboi's perfectness result (Lemma 3.3 in \cite{AC}) we use Lemma \ref{lem:tsuboi} of Tsuboi which is a more subtle perfectness result in the higher regularity. Finally, in the conclusion, we use Lemma \ref{lem:rootCk} instead of Lemma \ref{lem:roots}.
   
   \bigskip
   
   Now, we would like to extend the result to all groups $\Diff^k_0(M) \ (\Diff^k_0(M,\partial M))$ where $M$ is a compact smooth manifold (with boundary) with $C^k$ topology. 
   
   \medskip
   
   \begin{thm} \label{thm:diff(M)} If $M$ is a compact smooth manifold then the group $\Diff^k_0(M)$ is topologically finitely generated. Similarly, if $M$ is a compact smooth manifold with boundary then the topological rank of $\Diff^k_0(M,\partial M)$ is finite. 
  \end{thm}
  
   \medskip
   
   {\bf Proof.} We may and will assume that $M$ is connected. 
   
   First, we will consider {\em the case of a closed manifold}. Let $d = \mathrm{dim}M$, $\Omega _1, \dots , \Omega _m$ be open charts of $M$ covering the entire manifold $M$ (the open sets $\Omega _s, 1\leq s\leq m$ are diffeomorphic to a unit open ball in $\mathbb{R}^d$). In each $\Omega _s, 1\leq s\leq m$ we take a sequence of open sets $D^{(s)}_1, D^{(s)}_2, D^{(s)}_3, \dots $ accumulating to some point $p_s\in \Omega _s\backslash  (\displaystyle \mathop{\cup }_{n\geq 1}\overline {D^{(s)}_n})$ such that the closures $\overline {D^{(s)}_1}, \overline {D^{(s)}_2}, \overline {D^{(s)}_3}, \dots $ are mutually disjoint. Since $\Omega _s$ can be identified with an open unit ball, we can and will take the open sets $D^{(s)}_1, D^{(s)}_2, \dots $ to be open balls in $\Omega _s$ with centers lying on a segment and accumulating to $p$ from one side. Let $\phi _s, \psi _s , 1\leq s\leq m$ be diffeomorphisms of $M$  such that the conditions analogous to conditions (a-i)-(a-iii) of \cite{AC} hold, more precisely,
   
   \medskip

   (i) $Int(\mathrm{supp} (\psi _s)) = \Omega _s$, 
   
    \medskip
    
    (ii) $\phi _s^{-1}(D^{(s)}_n) = D^{(s)}_{n+1}, \forall n\geq 1$,
    
    \medskip
    
    (iii) the sets  $A_k = \psi ^{-k}D^{(s)}_1$ form an increasing chain of open sets in $D^{(s)}_1$ (i.e. $A_1\subseteq A_2 \subseteq \dots \subseteq D^{(s)}_1$) and $\displaystyle \mathop {\cup }_{k\geq 1}A_k = \Omega _s$.

  \medskip
  
  Notice that the group $\Diff^k_0(M)$ is generated by subgroups $G_s = \{\eta \in \mathrm{Diff}^k_0(M) \ | \ \mathrm{supp } (\eta ) \cup \mathrm{Im} (\eta ) \subseteq \overline{\Omega _s}\}$. Hence it suffices to build maps $f_s, g_s, u_s, v_s, h_s \in G_s$ such that $f_s, g_s, u_s, v_s, h_s, \phi _s, \psi _s$ generate $G_s$.  
  
  \medskip
  
  Now, fix $s\in \{1,\dots , m\}$ and let $\eta _1, \eta _2, \dots $ be dense sequence in the subset $\{[\omega _1\omega _2,\omega _3\omega _4] : \omega _i\in G_s, 1\leq i\leq 4\}$. We will identify the set $\Omega _s$ with the open unit ball in $\mathbb{R}^d$.  
  
  \medskip
  
  We define $f, g, u, v$ and $h$ on $M\backslash \displaystyle \mathop {\sqcup }_{n\geq 1}D^{(s)}_n$ to be an identity map and build these maps on the open sets $D^{(s)}_1, D^{(s)}_2, \dots $ inductively as in the proof of Theorem 2.1 from \cite{AC}. More precisely, let $F_0$ be a free group formally generated by $f,g,u,v,h$ and $(\epsilon _n)_{n\geq 1}$ be a positive sequence decreasing to zero.
  
  \medskip
  
   In the domain $D^{(s)}_1$ we define $f, g, u, v$ and $h$ such that all of them map $D^{(s)}_1$ to itself and all are $C^k$-flat at the boundary $\partial (D^{(s)}_1)$; moreover, for some word $X_1(f,g,u,v)$ in $F_0$ we have $\phi ^{-1}\psi ^{m_1}\phi X_1(f,g,u,v)\phi ^{-1}\psi ^{-m_1}\phi $ is $\epsilon _1$ close to $\eta _1$ for some sufficiently big $m_1$. 
  
  \medskip
  
  Suppose now $f, g, u, v$ and $h$ are defined on $\displaystyle \mathop{\sqcup }_{1\leq i\leq n}D^{(s)}_i$, as well as the words $X_i(f,g,u,v,h), Y_i(f,g), Z_i(u,v), 1\leq i\leq n$ are defined such that all the maps $f, g, u, v, h$ are mapping each of $D^{(s)}_i, 1\leq i\leq n$ to itself, and the following conditions analogous to (b-i)-(b-iv) of \cite{AC} hold:
  
  \medskip
  
  (i) for all $x\in \partial D^{(s)}_{i}$ and for all $l, j\in \{1, \dots , d\}$, $$(\partial f_l/\partial x_j)(x) = (\partial g_l/\partial x_j)(x) = (\partial u_l/\partial x_j)(x) = (\partial v_l/\partial x_j)(x) = (\partial h_l/\partial x_j)(x) = \delta _l^j$$ and for all $2\leq r\leq k, j_1, \dots, j_r, l \in \{1, \dots , d\}$, $$(\partial ^r\theta _l/\partial x_{j_1}\dots \partial x_{j_r})(x) = 0$$ where $\theta _l$ is any of the functions $f_l, g_l, u_l, v_l, h_l$; 
  
  \medskip
  
  (ii) $f, g, u, v$ and $h$ are $\epsilon _{i}$-close to the identity in $C^k(D^{(s)}_{i})$; 
  
  \medskip
  
  (iii) for some natural number $p_{i}$ we have  $$(\mathrm{supp}(h^{p}Y_j(f,g)h^{-p})\cap \mathrm{supp}(Z_j(u,v)))\cap D^{(s)}_{i} = \emptyset $$ for all $1\leq i\leq n, 1\leq j < i$ and $p\geq p_i$;
   
   \medskip
   
   (iv) for the word $X_{i} = [h^{p_{i}}Y_{i}(f,g)h^{p_{i}}, Z_{i}(u,v)]$ in $F_0$ we have $$\phi ^{-i}\psi ^{m_i}\phi ^{i}X_{i}(f,g,u,v,h)\phi ^{-i}\psi ^{-m_i}\phi ^{i}$$ is $\epsilon _{i}$-close to $\eta _{i}$ for some sufficiently big $m_i$.

  \medskip
  
  Then we define the diffeomorphisms $f = (f_1, \dots , f_d), g = (g_1, \dots , g_d), u = (u_1, \dots , u_d), v = (v_1, \dots , v_d), h = (h_1, \dots , h_d)$ on the domain $D^{(s)}_{n+1}$ such that all these maps are mapping $D^{(s)}_{n+1}$ to itself, and
  
  (i) for all $x\in \partial D^{(s)}_{n+1}$ and for all $l, j\in \{1, \dots , d\}$,  $$(\partial f_l/\partial x_j)(x) = (\partial g_l/\partial x_j)(x) = (\partial u_l/\partial x_j)(x) = (\partial v_l/\partial x_j)(x) = (\partial h_l/\partial x_j)(x) = \delta _l^j$$ and for all $2\leq r\leq k, j_1, \dots, j_r, l \in \{1, \dots , d\}$, $$(\partial ^r\theta _l/\partial x_{j_1}\dots \partial x_{j_r})(x) = 0$$ where $\theta _l$ is any of the functions $f_l, g_l, u_l, v_l, h_l$; 
  
  \medskip
  
  (ii) $f, g, u, v$ and $h$ are $\epsilon _{n+1}$-close to the identity in $C^k(D^{(s)}_{n+1})$;
  
   \medskip
   
   (iii) for some natural number $p_{n+1}$ we have  $$(\mathrm{supp}(h^{p}Y_i(f,g)h^{-p})\cap \mathrm{supp}(Z_i(u,v)))\cap D^{(s)}_{i} = \emptyset $$ for all $1\leq i\leq n, p\geq p_{n+1}$;
   
   \medskip
   
   (iv) for some words $Y_{n+1}(f,g), Z_{n+1}(u,v)$ and for the word $X_{n+1} = [h^{p_{n+1}}Y_{n+1}(f,g)h^{p_{n+1}}, Z_{n+1}(u,v)]$ in $F_0$ we have $$\phi ^{-(n+1)}\psi ^{m_{n+1}}\phi ^{n+1}X_{n+1}(f,g,u,v,h)\phi ^{-(n+1)}\psi ^{-m_{n+1}}\phi ^{n+1}$$ is $\epsilon _{n+1}$-close to $\eta _{n+1}$ for some sufficiently big $m_{n+1}$.
   
  \medskip  
  
  By induction, we extend $f, g, u, v$ and $h$ to the whole domain $\Omega _s$. By construction, the group generated by $\phi , \psi , f, g, u, v$ and $h$ is dense in $G_s$. Indeed, we can take the words $Y_{n}(f,g) = f^{M_1}g^{M_2}$ and $Z_{n}(u,v) = u^{M_3}v^{M_4}$ for sufficiently big $M_1, M_2, M_3, M_4$, and the claim now follows from Lemma \ref{lem:rootextlem} and Lemma \ref{lem:hrt}. 
  
  \medskip
  
  For the {\em case of a compact manifold with boundary,} let us first notice that the subgroup $\Diff ^k(M, \partial M)_{flat}$ of $C^k$-flat diffeomorphisms of $\Diff ^k(M, \partial M)$ is also topologically finitely generated by the treatment of the case of closed manifolds. Secondly, if $N$ is a closed manifold then the group
  
   $G_N = \{f\in \Diff _0^k (N\times I) \ | \ \forall x\in N, f(x,0) = (x,0), \ \mathrm{and} $ \ $ \ \forall x\in N, \forall t\in I, f(x,t) = (y,t), \ \mathrm{for \ some} \ y\in N\}$ is topologically finitely generated as well similarly to the case of a closed manifold. Now, by Collar Neighborhood Theorem, the boundary $\partial M$ has an open neighborhood in the form $\partial M\times [0,2)$ in $M$. Denote $N = \partial M$. Then both of the groups $G_N$ and $\Diff _c^k(M, \partial M)_{flat}$ have finite topological ranks. But these two groups generate $\Diff _0^k(M, \partial M)$ thus we obtain that the topological rank of the latter is finite. $\square $ 
  
  \medskip
  
  \begin{rem} It follows from the proof that the group Diff$^k_{+}(I)$ is generated by $3k + 7$ elements, and the group Diff$^k_0(M)$ is generated by $(3k+7)m$ elements where $m$ is the minimal number of charts covering $M$. With more careful analysis it is possible to decrease this bound significantly. Indeed, it is quite easy to make the bound
 independent of $m$; in the case of a closed manifold, one simply needs to take a diffeomorphism $\pi :M\to M$ and an open chart $\Omega _1$ such that $(\pi ^n(\Omega _1))_{n\geq 0}$ is an open cover of $M$ (similar arrangements can be made for manifolds with boundary as well). It is also possible to make the bound independent of $k$. We suspect that in the case of a closed manifold the minimal of generators can be made equal to two. It also follows from the construction that the generators can be chosen from an arbitrary small neighborhood of identity (thus the so called {\em local rank} is also finite), moreover, all the generators can be chosen from the $C^{\infty }$ class. 
  \end{rem}

 \vspace{1cm}
 
 \section{Roots of diffeomorphisms and perfectness results in diffeomorphism groups}
 
  In this section we will prove the Lemma \ref{lem:introroot}, and then generalize this result to a higher dimensional case. We will need the following
  
  \medskip
  
  \begin{defn} Let $l\geq 1$. A sequence $(t_i)_{1\leq i\leq N}$ of real numbers is called {\em $l$-quasi-monotone} if there exists $1 < i_1 \leq \dots \leq i_s < N$ for some $s\in \{1, \dots ,l\}$ such that each of the subsequences $(t_{i_j}, \dots , t_{i_{j+1}})$ is monotone where $0\leq j\leq s$ and $i_0 = 1, i_{s+1} = N$. A 1-quasi-monotone sequence will be simply called {\em quasi-monotone}.
  \end{defn}
  
  \medskip  
  
  In \cite{AC}, we prove the following lemma (we recommend that the reader studies the proof of it, since in this section, we are making a reference to the proof.)
 
 \begin{lem}\label{lem:roots} Let $f\in \mathrm{Diff}_{+}(I), \ f'(0) = f'(1) = 1$ and $r > 0$. Let also $f$ has no fixed point in $(0,1)$. Then there exists a natural number $N\geq 1$ and $g\in \mathrm{Diff}_{+}(I)$ such that $||g||_1 < r$ and $||g^N-f||_1 < r$. 
 \end{lem}
 
 \bigskip
 
  Now, we would like to extend the proof of the above lemma to the case when $f\in \mathrm{Diff}_{+}^k(I)$ hence we would like to find an approximate root $g$ of $f$ from a small neighborhood of identity in the $C^k$ metric. Here, we need estimations at all levels up to the $k$-th derivative. Thus, instead of inequalities {\bf (1)} and {\bf (2)} from the proof of Lemma 1.2 in \cite{AC}), we will need $k+1$ inequalities. As before, we will define the values of $g$ at the selected points $z_i, 0\leq i\leq nN$, then we will define the values of the derivatives at those points, then the values of the second derivatives and so on. Then, we extend $g$ to the whole interval $[0,1]$ such that $g'(x)$ will be monotone at each subinterval $(z_i, z_{i+1})$, moreover, the higher derivatives $g^{(j)}(z)$ do not deviate much from the values $g^{(j)}(z_i), g^{(j)}(z_{i+1})$ at the end points; let us also emphasize that we do not demand the monotonicity of the higher derivatives at the subintervals $(z_i, z_{i+1})$. In the proof of Lemma \ref{lem:roots}, we also used the $l$-quasi-monotonicity of the sequence $g'(z_0), \dots , g'(z_{nN})$, i.e. it can divided into at most $l$ (not necessarily strictly) monotone subsequences (in the simplest case of 1-quasi-monotonicity, this means that for all $1\leq i\leq nN-N$ there exists  $p\in \{i, \dots , i+N-1\}$ such that both of the sequences $(g'(z_i), \dots , g'(z_{i+p}))$ and  $(g'(z_{i+p}), \dots , g'(z_{i+N}))$ are monotone). This follows directly from the Mean Value Theorem applied to all $j$-th order derivatives, $0\leq j\leq k$, of iterates of $f^s, 1\leq s\leq n$ where $n$ will be chosen in the proof. Also, we will not impose a quasi-monotonicity condition on any of the sequences $g^{(j)}(z_0), \dots , g^{(j)}(z_{nN})$ for any $j\geq 2$.  
  
  \bigskip
 
  The following two lemmas are very technical although in essence they are very simple claims. These lemmas can also serve as a guide to the proof of the root approximation lemma, and alert the reader about the strategy of the proof.  To state the first lemma, we need to introduce some terminology.
  
  \medskip
  
   Let $u\in \mathrm{Diff}_{+}(I)$, $n, N$ be positive integers, and $0 = z_0 < z_1 < \dots < z_{nN-1} < z_{nN} = 1$ such that $z_{i} = \frac{i}{N^2}, 0\leq i\leq N$, and $z_{i+1} = u(z_i), N\leq i\leq nN-1$. In this case we say $(z_i)_{0\leq i\leq nN}$ is an {\em $(n,N)$-associated sequence} of $u$. For every $j\geq 1, 0\leq i\leq nN-j$, inductively on $j$, we define the quantities $$\overline{u^{(j+1)}(z_i)} = \frac{\overline{u^{(j)}(z_{i+1})} - \overline{u^{(j)}(z_i)}}{z_{i+1}-z_i}$$
where we assume that $\overline{u^{(0)}(z_i)} = z_{i+1}$ for all $0\leq i\leq nN$.  

  \medskip
  
  \begin{lem}\label{lem:buildg} Let $u\in \mathrm{Diff}_{+}(I)$ be a $C^{k+1}$-flat polynomial diffeomorphism, $k, n$ be positive integers. Then there exists $K > 0$ such that for all sufficiently big $N \geq 1$, and for all $\epsilon : \mathbb{N}\to (0,1)$ such that $\epsilon (m) \searrow 0$ as $m\to \infty $, if $(z_i)_{0\leq i\leq nN}$ is an  $(n,N)$-associated sequence of $u$ then there exists a diffeomorphism $g\in \Diff _{+}^{k+1}(I)$ such that $g(z_i) = z_{i+1}, 0\leq i\leq nN-1$, moreover, for all $0\leq i\leq nN-k-1, 2\leq j\leq k+1$, the following conditions hold:
  
  (i) $g'(x)$ is monotone on $(z_i, z_{i+1})$, 
  
  (ii) $|g'(z_i) - 1| < \frac{K}{N}$,
  
  (iii) $|g'(z_i) - \overline{u^{(1)}(z_i)}| < \epsilon (N)$,
  
  (iv) $|g^{(j)}(z_i)| < \frac{K}{N}$,
  
  (v) $|g^{(j)}(x)| < 2K\max \{g^{(j)}(z_i), g^{(j)}(z_{i+1})\}$ for all $x\in (z_i, z_{i+1})$.
\end{lem}

 \medskip
 
 {\bf Proof.} First, let us notice that, since $u$ is a $C^{k+1}$-flat polynomial, for some $K_0 > 0$ independent of $N$, $$|\overline{u^{(1)}(z_i)}-1| < \frac{K_0}{N}, 0\leq i\leq nN-1 \ \mathrm{and} \ |\overline{u^{(j)}(z_i)}| < \frac{K_0}{N}, 0\leq i\leq nN-j.$$
 
 Then condition (ii) follows from condition (iii). 
 
 Now, for all $0\leq i\leq nN-1$, let 
 
  \begin{displaymath}  g'(z_i) =  \left\{\begin{array}{lcr} \overline{u^{(1)}(z_i)} + \epsilon (N) &  \mathrm{if} \ \overline{u^{(1)}(z_i)} < \overline{u^{(1)}(z_{i+1})} \\ \overline{u^{(1)}(z_i)} - \epsilon (N) &  \mathrm{if} \ \overline{u^{(1)}(z_i)} > \overline{u^{(1)}(z_{i+1})} \\  \overline{u^{(1)}(z_i)} &  \mathrm{if} \ \overline{u^{(1)}(z_i)} = \overline{u^{(1)}(z_{i+1})}  \end{array}   \right. \end{displaymath} 
 
 \bigskip
 
 If $\epsilon (N)$ is sufficiently small then $g$ can be extended such that $g'(x)$ is monotone in the intervals $(z_i, z_{i+1})$. To satisfy conditions (iv) and (v), it suffices to satisfy $|g^{(k+1)}(x)| < \frac{K}{N}, 0\leq i\leq nN-k-1, x\in [0, 1]$, and this is straightforward. $\square $

  \medskip

  The second lemma is perhaps less obvious; we present a detailed proof of it. 

 \medskip
     
  \begin{lem}\label{lem:choice} Let $n$ be a positive integer, and $M > 1, r > 0$. Then there exists $K = K(M, n)$ such that for sufficiently big $N$, if 
  
   (i) $\frac{1}{M} < d_i^{(q)} < M, 1\leq i\leq nN, 1\leq j\leq N$,
   
   (ii) $|a_i| < M, 1\leq i\leq nN-n$,
   
   (iii) $|d_i^{q}-d_{i'}^{q'}| < \frac{M}{N}$, for all $1\leq i, i'\leq nN, 1\leq q, q'\leq N$ such that $\max\{|i-i'|, |q-q'|\} \leq 1$.
   
   (iv) $|a_i-a_{i+1}|< \frac{M}{N}, 1\leq i\leq nN-1$,
   
   then there exist $x_1, \dots , x_{nN}$ such that $|x_i| < \frac{K}{N}, 1\leq i\leq nN$ and $|\displaystyle \mathop{\sum }_{q=1}^Nd_{i}^{(q)}x_{i-1+q} - a_i| < r, 1\leq i\leq nN-n$. 
  \end{lem}
  
  {\bf Proof.} We may and will assume that $M > 2n$. Let us first observe that for arbitrary $x_1, \dots , x_{nN-1}$ we have  $$|\displaystyle \mathop{\sum }_{q=1}^Nd_{i}^{(q)}x_{i-1+q} - \displaystyle \mathop{\sum }_{q=1}^Nd_{i+1}^{(q)}x_{i+q}| < \displaystyle \mathop{\sum }_{q=2}^N|d_{i}^{(q)}-d_{i+1}^{(q)}||x_{i-1+q}| + |d_{i}^{(1)}x_{i}- d_{i+N}^{(N)}x_{i+N}| < $$ \ $$ <\displaystyle \mathop{\sum }_{q=2}^N\frac{M}{N}\frac{K_i}{N} + |d_{i}^{(1)}x_{i}- d_{i+N}^{(N)}x_{i+N}| < \frac{MK_i}{N} + |d_{i}^{(1)}x_{i}- d_{i+N}^{(N)}x_{i+N}| \ (\ast )$$ where $K_i = N\max \{|x_i|, \dots , |x_{i+N}|\}$.  
  
  \medskip
  
  Now, let $m = [M^4+1]$ and $N+1 = s_0 < s_1 < \dots < s_{m-1} < s_{m} = nN$ such that the quantities $s_{i+1}-s_{i}, 1\leq i\leq m-1$ differ from each other by at most one unit.
  
 \medskip   
  
  We will define the sequence $(x_p)_{1\leq i\leq nN}$ inductively on $p$ as follows: first, we choose $x_1 = x_2 = \dots = x_{N} \in (0, \frac{M^2}{N}]$ such that $\displaystyle \mathop{\sum }_{q=1}^Nd_{1}^{(q)}x_{q} = a_1$. 
  
  For each $p\in [s_i, s_{i+1}), 1\leq i\leq m-1$, we choose $x_p = \pm \frac{M^{4i}}{N}$ such that $x_p > 0$ if $\displaystyle \mathop{\sum }_{q=1}^Nd_{p-N}^{(q)}x_{p-N-1+q} < a_{p-N}$, and $x_p < 0$ otherwise.

  Thus to guarantee the inequality $|\displaystyle \mathop{\sum }_{q=1}^Nd_{p-N+1}^{(q)}x_{p-N+q} - a_{p-N+1}| < r $,
 at each step $p$, we make an adjustment (when defining $x_p$) depending if the previous sum is smaller or bigger than needed (i.e. whether or not $\displaystyle \mathop{\sum }_{q=1}^Nd_{p-N}^{(q)}x_{p-N-1+q} < a_{p-N}$). The size (absolute value) of the adjustment grows as the step $p$ grows, but it is not bigger than $\frac{M^{4m}}{N}$. Similar to the inequality $(\ast )$ one can show that at every step $p\in [s_i, s_{i+1})$ we have the inequality $$|\displaystyle \mathop{\sum }_{q=1}^{N}d_{p-N}^{(q)}x_{p-N-1+q}-\displaystyle \mathop{\sum }_{q=1}^{N-1}d_{p-N+1}^{(q)}x_{p-N+q}| < \frac{M^{4i-2}}{N} \ (\ast _2).$$
 
 Indeed, we have $$|\displaystyle \mathop{\sum }_{q=1}^{N}d_{p-N}^{(q)}x_{p-N-1+q}-\displaystyle \mathop{\sum }_{q=1}^{N-1}d_{p-N+1}^{(q)}x_{p-N+q}| \leq d_{p-N}^{(1)}|x_{p-N}| + \displaystyle \mathop{\sum }_{q=1}^{N-1}|d_{p-N}^{(q)}-d_{p-N+1}^{(q)}| |x_{p-N+q}|$$ \ $$\leq d_{p-N}^{(1)}|x_{p-N}| + \displaystyle \mathop{\sum }_{p-N+1+q\leq N}|d_{p-N}^{(q)}-d_{p-N+1}^{(q)}| |x_{p-N+q}| + $$ \ $$\displaystyle \mathop{\sum } _{i+1}^t\displaystyle \mathop{\sum }_{p-N+1+q\in [s_i, s_{i+1})}|d_{p-N}^{(q)}-d_{p-N+1}^{(q)}| |x_{p-N+q}| \leq$$ \   $$M\frac{M^2}{N} + \frac{nN}{M^4}\frac{M}{N}(\frac{M^4}{N} + \frac{M^{8}}{N} + \dots + \frac{M^{4i}}{N}) = \frac{M^3}{N} + \frac{n}{M^3N} \frac{M^{4(i+1)}-1}{M^4-1}$$ \ $$\leq \frac{M^3}{N} + \frac{n}{M^3N} M^{4i} \leq \frac{M^{4i-2}}{N}.$$
 
  Hence $$|\displaystyle \mathop{\sum }_{q=1}^{N}d_{p-N}^{(q)}x_{p-N-1+q}-\displaystyle \mathop{\sum }_{q=1}^{N}d_{p-N+1}^{(q)}x_{p-N+q}| \geq \frac{M^{4i-1}}{N} - \frac{M^{4i-2}}{N} > \frac{M^{4i-2}}{N}.$$ 
 
  Thus we can guarantee the inequality $|\displaystyle \mathop{\sum }_{q=1}^Nd_{i}^{(q)}x_{i+q} - a_{i}| < r$ for all $1\leq i\leq nN-N$ by satisfying the inequality $|x_1| \leq |x_2| \leq \dots \leq |x_{nN-1}| \leq \frac{M^{4m}}{N}$. Thus it suffices to take $K = M^{4M^4+4}$. $\square $.
 
 \medskip
 
 Now we are ready to prove the main lemma about the approximation of roots in the diffeomorphism group $\mathrm{Diff}_{+}^k(I)$.  
   
   \begin{lem}\label{lem:rootCk} Let $f\in \mathrm{Diff}_{+}^k(I), \ f'(0) = f'(1) = 1, f^{(j)}(0) = f^{(j)}(1) = 0, 2\leq j\leq k$ and $r > 0$. Let also $f$ has no fixed point in (0,1). Then there exists $N\geq 1$ and $g\in \mathrm{Diff}_{+}^k(I)$ such that $||g||_k < r$ and $||g^N-f||_k < r$. 
 \end{lem}
 
 \medskip
 
 {\bf Proof.} By replacing $f$ with a $C^{k+1}$-flat polynomial diffeomorphism from its $r/4$ neighborhood, we may assume that $f\in \mathrm{Diff}_{+}^{k+1}(I)$, and the functions  $f'(x)-1$ and $f^{(j)}(x), 2\leq j\leq k$ change their sign at most $l$ times. Let $C = 2 + ||f(x)-x||_{k+1}$, $n$ be sufficiently big with respect to $l$, and $0 < x_0 < x_1 < \dots < x_n < 1$ such that $$\max \{x_0, 1-x_n\} < r/2, f(x_i) = x_{i+1}, 0\leq i\leq n-1,$$ \ moreover, we have  $$\displaystyle \mathop{\sup }_{x\in [0, x_1]\cup [x_{n-1},1]}|f'(x) - 1| < r/2 \ \mathrm{and} \ \displaystyle \mathop{\sup }_{x\in [0, x_1]\cup [x_{n-1},1], 2\leq j\leq k+1}|f^{(j)}(x)| < r/2.$$  
  
 \medskip
 
 Let $N\geq 1$ be sufficiently big with respect to $\max\{n, C, l\}$, and $z_0, z_1, z_2,$ \ $\dots , z_{nN}$ be a sequence such that $z_{iN} = x_i, 0\leq i\leq n, z_{k} = x_0 + \frac{x_1-x_0}{N}k, 0\leq k\leq N$, and $z_{iN+j} = f(z_{iN-N+j}), 1\leq i\leq n-1, 1\leq j\leq N$. The subsequence $(z_1, z_2, \dots , z_{N-1})$ will be called {\em the initial subsequence}; it determines the entire sequence $(z_i)_{0\leq i\leq nN}$. 
 
 \medskip
 
 We define the diffeomorphism $g\in \mathrm{Diff}_{+}(I)$ as follows. Firstly, we let $g(z_i) = z_{i+1}, 0\leq i\leq nN-1$.
 Let us observe that we already have $g^N(z_i) = z_{i+N} = f(z_i), 0\leq i\leq nN-N$. Notice that as $N\to \infty $ the quantity $\displaystyle \mathop{\max }_{0\leq i\leq N-1}\frac{z_{i+2}-z_{i+1}}{z_{i+1}-z{i}}$ tends to 1, by Mean Value Theorem. Hence, we have $||g^N(x) - f(x)||_0 < r$ {\bf (1)} (i.e. for all possible extensions of $g$).

  \medskip

 Now we need to define $(g^N)^{(j)}(z_i), \ 0\leq i\leq nN-N, 1\leq j\leq k$ such that for a suitable extension of $g$ we could still claim that $||g^N(x) - f(x)||_k < r$. (Obviously, the latter cannot hold for an arbitrary extension of $g$.) We also need to make sure that $g$ will be small (i.e. close to the identity) in the $C^k$ norm.
 
 \medskip
 
 By chain rule, to define the quantities $(g^N)^{(j)}(z_i), \ 0\leq i\leq nN-N, 1\leq j\leq k$, it suffices to define the values $g^{(j)}(z_i), \ 0\leq i\leq nN-N, 1\leq j\leq k$ of the derivatives of the function $g$ up to the $k$-th order. Let us observe that for all $0\leq i\leq nN-N$, we have $$(g^N)''(z_i) = g''(z_i)\frac{\Pi }{g'(z_i)} + g''(z_{i+1})\frac{\Pi g'(z_{i})}{g'(z_{i+1})} + g''(z_{i+2})\frac{\Pi g'(z_{i})g'(z_{i+1})}{g'(z_{i+2})} + $$ $$+ \dots + g''(z_{i+N-1})\frac {\Pi g'(z_{i})g'(z_{i+1})\dots g'(z_{i+N-2})}{g'(z_{i+N-1})} \ (\ast )$$
 
 where $\Pi = g'(z_i)g'(z_{i+1})\dots g'(z_{i+N-1})$.
 
 \medskip
 
 The right-hand side of the formula $(\ast )$ defines a polynomial $P_2$ (independent of $i$) of $2N$ variables \footnote{More precisely, the polynomial $P_2$ is given by the formula $$P_2(X_1, \dots , X_N, Y_1, \dots , Y_N) = R_1(X_1, \dots , X_N)Y_1 + \dots + R_N(X_1, \dots , X_N)Y_N$$ where $R_1, \dots , R_N$ are polynomials in $N$ variables $X_1, \dots , X_N$ defined as follows: $$R_1(X_1, \dots , X_N) = X_2X_3\dots X_N = \frac{\omega }{X_1},$$ \  $$R_2(X_1, \dots , X_N) = X_1^2X_3\dots X_N = \frac{\omega X_1}{X_2},$$ \ $$R_3(X_1, \dots , X_N) = X_1^2X_2^2X_4\dots X_N = \frac{\omega X_1X_2}{X_3}, \dots ,$$ \ $$R_N(X_1, \dots , X_N) = X_1^2X_2^2\dots X_{N-1}^2 = \frac{\omega X_1X_2\dots X_{N-1}}{X_N}.$$ Here, $\omega = X_1\dots X_N$.}
 
  such that $$(g^N)''(z_i) = P_2(g'(z_i), \dots , g'(z_{i+N-1}), g''(z_i), \dots , g''(z_{i+N-1})), 0\leq i\leq nN-N.$$

 \medskip
 
 Let us also observe that if $g$ is of class $C^j, j\geq 3$ then $$(g^N)^{(j)}(z_i) = g^{(j)}(z_i)\frac{\Pi }{g'(z_i)} + \dots +  g^{(j)}(z_{i+N-1})\frac {\Pi g'(z_{i})g'(z_{i+1})\dots g'(z_{i+N-2})}{g'(z_{i+N-1})} + A_1 + \dots + A_d \ (\ast _j)$$ where $d\leq N^{2j}$ for sufficiently big $N$, and each of the terms $A_p, 1\leq p\leq d$ consists of a product of at least $N-1$ and at most $N^j$ derivatives $g^{(j-l)}(z_{i+s}), 0 < l\leq j-1, 0\leq s\leq N-1$, moreover, at least two and at most $j$ of the derivatives in the product are values of derivatives of order $j-l\in \{2, \dots , j-1\}$. 
 
 \medskip
 
 Again, the right-hand side of the formula $(\ast _j)$ defines a polynomial $P_j$ of $jN$ variables such that $$(g^N)^{(j)}(z_i) = P_j(g'(z_i), \dots , g'(z_{i+N-1}), \dots , g^{(j)}(z_i), \dots , g^{(j)}(z_{i+N-1})), 0\leq i\leq nN-N.$$ We make an important observation that $$P_j(g'(z_i), \dots , g'(z_{i+N-1}), \dots , g^{(j)}(z_i), \dots , g^{(j)}(z_{i+N-1})) = $$ \ $$ = P_2(g'(z_i), \dots , g'(z_{i+N-1}), g^{(j)}(z_i), \dots , g^{(j)}(z_{i+N-1})) + A_1 + \dots + A_d \ (\ast )$$
 
 \medskip
 
 Now, we would like to introduce the quantities $\overline{g^{(j)}(z_i)}, 0\leq j\leq k+1, 0\leq i\leq nN-k$ by defining them inductively on $j$ as follows: $$\overline {g^{(j+1)}(z_i)} = \frac{\overline {g^{(j)}(z_{i+1})}- \overline {g^{(j)}(z_i)}}{z_{i+1}-z_i}$$ and $\overline {g^{(0)}(z_i)} = g(z_i)$.
 
 \medskip
 
 Then, using the formula $(\ast _j)$ and Lemma \ref{lem:choice}, inductively on $j$, for sufficiently big $N$, we can choose a constant $K = K(C, n, r, l) > C^{k+1}$, the values $g^{(j)}(z_i), 0\leq i\leq nN-N$ such that the following conditions hold 
 
 \medskip
 
 (c-i) $|\overline {g'(z_i)}-1| < \frac{K}{N}$,
 
 \medskip
 
 (c-ii) $|\overline {g^{(j)}(z_i)}| < \frac{K}{N}$ for all $2\leq j\leq k+1$,
 
 \medskip
 
 (c-iii) the sequence $\overline {g'(z_i)}, \dots , \overline {g'(z_{i+N})}$ is $l$-quasi-monotone;
 
 \medskip
 
 (c-iv) $\frac{1}{K}|\overline {g^{(j)}(z_i)}| < |g^{(j)}(z_i)| < K|\overline {g^{(j)}(z_i)}|$ for all $1\leq j\leq k+1$,
 
 \medskip
 
 (c-v) $|P_j(g'(z_i), \dots , g'(z_{i+N-1}), \dots , g^{(j)}(z_i), \dots , g^{(j)}(z_{i+N-1})) -f^{(j)}(z_i)| < \frac{r}{2}$
 
 \medskip
 
 (c-vi) $g'(z_i) = \frac{z_{i+2}-z_{i+1}}{z_{i+1}-z_i} + \epsilon _i^{(N)}, 0\leq i\leq nN-2$, where $\mathrm{sgn}(\epsilon _i^{(N)}) = \mathrm{sgn}(\overline{g'(z_i)} - \overline{g'(z_{i+1})})$.
 
 \medskip
 
 (c-vii) the sequences $(g'(z_i), \dots , g'(z_{i+N}))$ are $l$-quasi-monotone.

 \medskip
 
 Indeed, conditions (c-i) and (c-ii)  for the sequence $\overline {g^{(j)}(z_i)}, 1\leq i\leq nN$ for all $1\leq j\leq k+1$, follow from the fact that all of the functions $(f^{(j)})^s, 0\leq j\leq k, 1\leq s\leq n$ are polynomials (i.e. all $n$ iterates of the function $f^{(j)}$ for all $0\leq j\leq k$). By the choice of $f$ and by the Mean Value Theorem, the monotonicity of the sequence $(\overline{g'(z_i)}, \dots , \overline{g'(z_{i+N})})$ is violated in at most $l$ terms. 
 
 \medskip
 
  After this, for all $2\leq j\leq k$, let $$d_{i,j}^{(q)} = R_q(g'(x_i), \dots , g'(x_{i+N-1}), g^{(j)}(x_i), \dots ,  g^{(j)}(x_{i+N-1}).$$ Let also $$a_i = P_2(g'(x_i), \dots , g'(x_{i+N-1}), g^{(j)}(x_i), \dots ,  g^{(j)}(x_{i+N-1})),$$ i.e. we have $$a_i = P_j(g'(x_i), \dots , g'(x_{i+N-1}), g''(x_i), \dots , g''(x_{i+N-1}), \dots , g^{(j)}(x_i), \dots ,  g^{(j)}(x_{i+N-1})$$ \ $$ - A_1 - \dots - A_d.$$ Then for each fixed $j\in \{2, \dots , k\}$, the quantities $d_{i}^{(q)}: = d_{i,j}^{(q)}$ satisfy the conditions of Lemma \ref{lem:choice} for some positive constants $M$ and $K(M, l,r, n)$ where the latter is independent of $N$. By this lemma and by Lemma \ref{lem:buildg}, using the formulas $(\ast _j)$, inductively on $j$, we can guarantee the conditions (c-iv) and (c-v). By choosing $\epsilon ^{(N)}$ sufficiently small and defining $g'(z_i), 1\leq i\leq nN$ as in condition (c-vi), by the $l$-quasi-monotonicity of $\overline {g'(z_i)}, 1\leq i\leq nN$ we also obtain the $l$-quasi-monotonicity of $g'(z_i), 1\leq i\leq nN$, which is the condition (c-vii) (i.e. the monotonicity of the sequence $(\overline{g'(z_i)}, \dots , \overline{g'(z_{i+N})})$ is violated in at most $l$ terms). 
 
 \medskip
 
 For simplicity, let us assume that the sequences $(g'(z_i), \dots , g'(z_{i+N}))$ are quasi-monotone, i.e. for some $p\in \{i, \dots , i+N-1\}$, both of the subsequences $(g'(z_i), \dots , g'(z_{i+N}))$ and  $(g'(z_{i+p}), \dots , g'(z_{i+N}))$ are monotone (again, the general case is similar). 
 
 \medskip      
 
 Now we need to extend the definition of $g(x)$ to all $x\in [0,1]$ such that we have $$|g'(x) - 1| < r {\bf (2)}, |g^{(s)}(x)| < r \ {\bf (3)} \ \mathrm{and} \  |(g^N)^{(j)}(x) - f^{(j)}(x)| < r \ {\bf (4)}$$ for all $x\in [0,1]$ and $2\leq s\leq k, 1\leq j\leq k$. By condition (c-v), for sufficiently big $N$, and sufficiently small $\epsilon ^{(N)}: = \displaystyle \mathop{\max }_{1\leq i\leq N}|\epsilon _i^{(N)}|$, we can extend the function $g$ at each of the intervals $(z_i, z_{i+1})$ such that $g'(z)$ is monotone on $(z_i, z_{i+1})$. Also, by Lemma \ref{lem:buildg} and by Mean Value Theorem, the extension can also be made to satisfy  $$|g^{(j)}(z)| <  2K\max \{g^{(j)}(z_i), g^{(j)}(z_{i+1})\} \ {\bf (5)}$$ and 
 
 $$|g^{(j)}(z) - g^{(j)}(z_i)| < \frac{2K^2}{N}|z_{i+1}-z_i| {\bf (6)}$$  for all $2\leq j\leq k$.
 
 \medskip
 
 By monotonicity of $g'(z)$ on each subinterval $(z_i, z_{i+1})$ and by condition (c-i) we obtain the inequality {\bf (2)}. From the inequality {\bf (5)} and condition (c-ii), we obtain the inequality {\bf (3)}. 
 
 \medskip
 
 For the inequality {\bf (4)}, let us first point out that, by condition (c-iv), we have this inequality for each value $z = z_i, 0\leq i\leq nN-k$. To obtain it for all possible values it suffices to show that $$|(g^N)^{(j)}(z) - (g^N)^{(j)}(z_i)| < r/2 \ {\bf (7)}$$ for all $z\in (z_i, z_{i+1})$. 
 
 \medskip
 
 For simplicity, we will first prove the inequality {\bf (7)} for $j = 2$. (for $j = 1$, it is already proved in the proof of Lemma 3.2 in \cite{AC}). The right-hand side of the formula $(\ast )$ also defines monomials $R_1, \dots , R_N$ on $N$ variables such that we can write $$(g^N)''(z_i) = g''(z_i)R_1(g'(z_i), \dots ,g'(z_{i+N-1})) + g''(z_{i+1})R_2(g'(z_i), \dots ,g'(z_{i+N-1})) +$$ \ $$ + \dots + g''(z_{i+N-1})R_N(g'(z_i), \dots ,g'(z_{i+N-1})).$$ \footnote{Let us recall that we have $$R_1(u_1, \dots , u_N) = u_2u_3\dots u_N,  \ R_2(u_1, \dots , u_N) = u_1^2u_3u_4\dots u_N,$$ \ $$R_3(u_1, \dots , u_N) = u_1^2u_2^2u_4u_5\dots u_N$$ and so on.}
 
 For $z\in (z_i, z_{i+1})$, we can also write $$(g^N)''(z) = g''(t_i)R_1(g'(t_i), \dots ,g'(t_{i+N-1})) + g''(t_{i+1})R_2(g'(t_i), \dots ,g'(t_{i+N-1})) + $$ \ $$+ \dots + g''(t_{i+N-1})R_N(g'(t_i), \dots ,g'(t_{i+N-1}))$$ where $t_m\in (z_m, z_{m+1})$ for all $m\in \{i, i+1, \dots , i+N-1\}$. 
 
 \medskip
 
 Then, taking into account that $(1+\frac{K}{N})^N\nearrow e^K$ we obtain \  $$|(g^N)''(z_i)-(g^N)''(z)| \leq \displaystyle \mathop{\sum }_{1\leq m\leq N}R_{m+1}(g'(z_i), \dots ,g'(z_{i+N-1}))|g''(z_{i-1+m})-g''(t_{i-1+m})| +$$ \ $$ +\displaystyle \mathop{\sum }_{1\leq m\leq N}|R_{m+1}(g'(z_i), \dots ,g'(z_{i+N-1}))-R_{m+1}(g'(t_i), \dots ,g'(t_{i+N-1}))|\cdot |g''(t_{i-1+m})| \leq $$ \ $$ \leq \displaystyle \mathop{\sum }_{1\leq m\leq N}R_{m+1}(g'(z_i), \dots ,g'(z_{i+N-1}))|z_{i+m}-z_{i-1+m}|\frac{2K^2}{N} + $$ \ $$+ \displaystyle \mathop{\sum }_{1\leq m\leq N}|R_{m+1}(g'(z_i), \dots ,g'(z_{i+N-1}))-R_{m+1}(g'(t_i), \dots ,g'(t_{i+N-1}))|\frac{2K^2}{N}$$ \ $$\leq e^{2K}\frac{2K^3}{N} + \displaystyle \mathop{\sum }_{1\leq m\leq N}C^2|z_{i+m}-z_{i-1+m}|\frac{2K^2}{N}  \leq 4e^{2K}\frac{K^3}{N}$$ 
 
 \medskip
 
 Thus we proved the inequality {\bf (7)} for $j = 2$. The proof for a general $j$ is similar. Let us call a monomial {\em special} if it is of the form $R(u_1, \dots , u_N) = u_pu_{p+1}\dots u_q$ for some $1\leq p \leq q \leq N$. Then we can write $$(g^N)^{(j)}(z_i) = g^{(j)}(z_i)R_1(g'(z_i), \dots ,g'(z_{i+N-1})) + \dots + $$\ $$ g^{(j)}(z_{i+N-1})R_N(g'(z_i), \dots ,g'(z_{i+N-1})) + \displaystyle \mathop{\sum }_{\alpha \in S}B_{\alpha  }R_{\alpha }(g'(z_i), \dots ,g'(z_{i+N-1}))$$ where $|S| = d$, and the terms $B_{\alpha }R_{\alpha }(g'(z_i), \dots ,g'(z_{i+N-1})) = A_{\alpha }, 1\leq \alpha \leq d$ are the terms $A_1, \dots , A_d$ from the formula $\ast )$. Thus each $B_{\alpha }$ is a product of at least two and at most $j$ derivatives $g^{(j-l)}(z_{i+s})$ where $1\leq l\leq j-2, 0\leq s\leq N-1$, and $R_{\alpha }$ is a product of at most $2j$ special monomials. Furthermore, the terms have the following finer properties: for each $2\leq q\leq j$, let $C_q$ denotes the set of ordered $q$-tuples $\omega = (k_1, \dots , k_q)$ (with repetitions allowed)  of the set $\{i, i+1, \dots , i+N-1\}$. Then one can partition $S = \sqcup _{1\leq q\leq j} D_q$ such that $|D_q| \leq (jN)^q|C_q|$ with a surjection $\xi _q:D_q\to C_q, 1\leq q\leq j$ where if $\alpha \in D_q$ and $\xi (\alpha ) = (k_1, \dots , k_q)$ then the term $B_{\alpha }$ consists of the product $g^{(j_1)}(x_{i+k_1})\dots g^{(j_q)}(x_{i+k_q})$ for some $j_1, \dots , j_q\in \{2, \dots , j\}$, moreover,  for all $x\in C_q$ we have $|\xi _q^{-1}(x)| \leq (jN)^q$. We let $$a_i = \displaystyle \mathop{\sum }_{\alpha \in S}A_{\alpha }R_{\alpha}(g'(z_i), \dots ,g'(z_{i+N-1}))$$ for all $0\leq i\leq nN-n$. By properties (c-i), (c-ii) and (c-iv), the quantities $d_{i,j}^{(q)}, 1\leq i\leq nN-n, 2\leq j\leq k, 1\leq q\leq N$ and $a_i, 0\leq i\leq nN-n$ satisfy the conditions of Lemma \ref{lem:choice}.  Then, by estimating the difference $|(g^N)^{(j)}(z_i) - (g^N)^{(j)}(z)|$ we can derive the inequality similarly to the case of $j = 2$  $\square $
 
 \bigskip
 
 Now we would like to quote a perfectness result from \cite{T} which we have used in Section 2. 
 
 \begin{lem}\label{lem:tsuboi} For all $f\in \mathrm{Diff}_{c}^{\infty }(I)$ there exists $u_1, \dots , u_8\in \mathrm{Diff}_{c}^{\infty }(I)$ such that $f = [u_1, u_2][u_3,u_4][u_5,u_6][u_7,u_8]$.
 \end{lem}
 
 \medskip
 
 {\bf Proof.} See \cite{T}.  
 
 \medskip
 
 Lemma \ref{lem:rootCk} and Lemma \ref{lem:tsuboi} are the results we needed for the proof of Theorem \ref{thm:main}. Now we need to extend these results to higher dimensions. Let $d\geq 2$, $\Omega $ be an open unit ball in $\mathbb{R}^d$, and $I^d = (0,1)^d$ be an open unit cube in $\mathbb{R}^d$. We would like to introduce a key notion of an {\em elementary diffeomorphism}.  
 
 \medskip
 
 \begin{defn} A diffeomorphism  $f\in \mathrm{Diff}^{\infty }_{c}(\Omega )$ is called an elementary diffeomorphism if there exists $\phi :\overline{I^d}$ to $\overline{\Omega }$ such that $\phi $ is a restriction of a $C^{\infty }$ diffeomorphism  $\Phi :\mathbb{R}^d\to \mathbb{R}^d$, and for all $(x,t)\in I^{d-1}\times I$ we have $\phi  ^{-1}f\phi  (x,t) = (x, \psi _x(t))$ where for all $x\in I^{d-1}$, the support of the diffeomorphism $\psi _x\in C_c^{\infty }(I)$ is connected. More generally, for $k\in \{1, \dots , d-1\}$, we say $f$ is a $k$-elementary diffeomorphism if there exists $\phi :\overline{I^d}\to \overline{\Omega }$ as a restriction of a $C^{\infty }$ diffeomorphism $\Phi :\mathbb{R}^d\to \mathbb{R}^d$ such that for all $(x,u)\in I^{d-k}\times I^k$ we have $\phi ^{-1}f\phi _(x,u) = (x, \psi _x(u))$ where, for all $x\in I^{d-k}$, the support of the diffeomorphism $\psi _x\in C_c^{\infty }(I)$ is connected.   
 \end{defn}

 \medskip
 
 \begin{prop} \label{lem:approx} Let $\Omega $ be an open unit ball in $\mathbb{R}^d$, $r > 0$ and $g\in \mathrm{Diff}^k_{c}(\Omega )$ where $k$ is a positive integer. Then there exists a diffeomorphism $h\in \mathrm{Diff}^k_{c}(\Omega )$ such that $h\in B_r(g)$ and $h$ is a product of elementary diffeomorphisms.
 \end{prop}
 
  {\bf Proof.} Let $S$ be unit sphere i.e. the boundary of $\Omega $, and $\tilde{S}$ be the quotient of $S$ obtained by identifying the antipodal points. For every unit vector ${\bf x}\in \mathbb{R}^d$, we will write $\mathcal{F}_{{\bf x}}$ to denote the foliation of $\Omega $ with straight segments parallel to ${\bf x}$. A diffeomorphism $f\in \mathrm{Diff}^k_{c}(\Omega )$ maps the foliation $\mathcal{F}_{{\bf x}} = \{L_z\}_{z\in \tilde{S}}$ to a foliation $f\mathcal{F}_{{\bf x}} = \{f(L_z)\}_{z\in \tilde{S}}$ where for all $z\in \tilde{S}$, the leaf $f(L_z)$ shares the same ends with $L_z$. 
 
 \medskip
 
 We will fix a unit vector ${\bf x_0}$. For simplicity, let us assume that $d = 2$. We will also assume that the interior of the support of $f$ is diffeomorphic to an open ball. Then the leaves $L_z$ and $f(L_z)$ bound a unique domain in $\Omega $ which we will denote by $D_z$. 
 
 \medskip
 
 An arc $C = (p,q)$ of a leave $f(L)$ of the foliation $f\mathcal{F}_{{\bf x_0}}$ will be called {\em maximal} if $C$ is {\em convex}, i.e. any segment connecting any two points on $C$ is totally inside the domain bounded by $L$ and $f(L)$, moreover, $C$ is not properly contained in any other convex arc.
 
 \medskip
 
 Let $\mathcal{A}_f$ denotes the set of all maximal arcs of $f\mathcal{F}_{{\bf x_0}}$. An element $a\in \mathcal{A}_f$ belongs to a certain leaf $L_z$, and will also write $D(a)$ to denote the domain $D_z$. We will introduce an equivalence relation on $\mathcal{A}_f$ as follows: the elements $a_1, a_2\in \mathcal{A}_f$ are called equivalent if there exists a continuous family of maximal arcs starting at $a_1$ and ending at $a_2$. The set of equivalence classes of $\mathcal{A}_f$ will be denoted with $A_f$. For any $a\in \mathcal{A}_f$ the union of all maximal arcs which can be reached from $a$ by a continuous family of maximal arcs will be denoted by $\Omega (a)$.  
 
 \medskip
 
 By assumption, the support of $f$ is diffeomorphic to a ball. This allows us to introduce a partial order in the set $A_f$. Let $\alpha _1, \alpha _2\in A_f$. We say $\alpha _1$ {\em dominates} $\alpha _2$ if for any $\epsilon > 0$ there exist representatives $a_1, a_2\in \mathcal{A}_f$ of $\alpha _1, \alpha _2$ respectively such that the arcs $a_1, a_2$ are $\epsilon $-close to each other, and $a_2$ belongs to the domain $D(a_1)$. Then, we say $\alpha _2 \prec \alpha _1$ if there exists a sequence $(\beta _1, \dots , \beta _n)$ such that $\beta _{i+1}$ dominates $\beta _i$ for all $0\leq i\leq n-1$, and $\beta _0 = \alpha _2, \beta _n = \alpha _1$.         
 
 \medskip
 
 Let us notice that if a diffeomorphism $f\in \mathrm{Diff}^k_{c}(\Omega )$ acts along the leaves of $\mathcal{F}_{{\bf x}}$ then it is an elementary diffeomorphism provided that the interior of the support is diffeomorphic to a ball.
 
 \medskip
 
 Now we choose $f$ from the neighborhood $B_{r/2}(g)$ in the space $\mathrm{Diff}^k_{c}(\Omega )$ such that $\mathrm{supp}f$ is diffeomorphic to a ball and the restriction of $f$ to this support coincides with a polynomial diffeomorphism. Then the set $A_f$ is finite. Let $n = |A_f|$ and $\alpha $ be a maximal element of $A_f$. Then there exists an elementary diffeomorphism $h_1$ such that $\mathrm{supp}(h_1)\subseteq \Omega (\alpha )$ and $|A_{h_1f}| < n$. Continuing the process inductively we find elementary diffeomorphisms $h_1, \dots , h_{n-1}$ such that the diffeomorphism $h_1\dots h_{n-1}f$ has only one equivalence class of maximal arcs, i.e. $|A_{h_1\dots h_{n-1}f}| = 1$. Then it is elementary. Thus $f$ is a product of elementary diffeomorphisms.   
 
 \medskip
 
 The case of higher dimensions $d > 2$ is similar. We replace the concept of maximal arc with a concept of maximal $k$-discs. Then we first show that $g$ can be approximated by a product of finitely many $(d-1)$-elementary diffemorphisms. Then we can approximate each of these $(d-1)$-elementary diffemorphisms with a product of finitely many $(d-2)$-elementary diffemorphisms, and proceed by induction.  $\square $   
 
 \bigskip
 
 The following lemma follows from the proof of Lemma \ref{lem:rootCk} as it is straightforward to extend the construction there (for the interval $I$) smoothly in the vertical direction of the cube $I^d = I\times I^{d-1}$.
 
 \medskip
 
 \begin{lem}\label{lem:rootextlem} Let $\Omega $ be an open unit ball in $\mathbb{R}^d$, $r > 0$ and $F\in \mathrm{Diff}^k_{c}(\Omega )$ be an elementary diffeomorphism. Then there exists $N \geq 1$ and $g\in \mathrm{Diff}^k_{c}(\Omega )$ such that $||g||_k < r$ and $||g^N-F||_k < r$.
 \end{lem} 
 
 \medskip
 
 {\bf Proof.} Indeed, in the proof of Lemma \ref{lem:rootCk},  By replacing $F$ with another diffeomorphism from its $r/4$-neghborhood if necessary, we may assume that there exists $N_0$ such that for all $x\in I^{d-1}$ and for all diffeomorphisms $F_x(t) = F(x,t)$ of the interval $I$ the number $N$ from the proof of Lemma \ref{lem:rootCk} can be taken equal to $N_0$. We will fix this $N_0$ for the rest of the proof.
   
   \medskip
   
   We also observe that the construction of the $N_0$-th root can be made continuous in $f$, i.e. for all $\epsilon > 0$ there exists $\delta = \delta (\epsilon ) > 0$ such that if $f_1$ is a $C^k$-flat diffeomorphism of $I$ and $||f_1-f||_k<\delta $, then one can construct a diffeomorphism $g_1$ such that $||g_1 - g||_k < \epsilon $ and $||g_1^{N_0}-f_1||_k < r/4$. 
   
   \medskip
   
   Now, let $\epsilon _0 > 0$ be such that if $||g-g_1||_k < \epsilon _0$ then $||g^{N_0}-g_1^{N_0}||_k < r/4$. Let also $\delta _0 = \delta (\epsilon _0)$. We can and will also assume that $\epsilon _0 < r/4$.

 \medskip
 
 Let $\delta _1 > 0$ be such that if $x, y\in I^{d-1}$ are at a distance at most $\delta _1$ apart, in the $l^1$-norm $||.||$ of $\mathbb{R}^{d-1}$, then $||F_x-F_y||_k < \delta _0$. Let $S = \frac{1}{m}\mathbb{Z}^{d-1}\cap I^{d-1}$ where $m$ is big enough that $S$ is a $\frac{1}{2}\delta _1$-net in $I^{d-1}$ (in the $l^1$-norm $||.||$ of $\mathbb{R}^{d-1}$). For all $x\in S$, we can choose a diffeomorphism $g_x(t)$ from an  $\epsilon _0$-neighborhood of the identity in the space $C^k(I)$ such that $||g_x^{N_0}-F_x||_k < r/4$; moreover, for all $x,y\in S$, if the distance $||x-y||$ is less than $\delta _1$ then $||g_x-g_y|| < \epsilon _0$. 
 
 \medskip
 
 We need to define the diffeomorphisms $g_x$ for all $x\in I^{d-1}$ (not just for $x\in S$). But for all $x,y\in S$, if the distance $||x-y||$ is less than $\delta _1$ then, in the space $C^k(I)$, the diffeomorphisms $g_x$ and $g_y$ can be connected by a path in the ball of radius $\epsilon _0$ centered at the origin. Thus we can extend a finite collection $\{g_x\}_{x\in S}$ to a diffeomorphism $g\in \mathrm{Diff}^k_{c}(I\times I^{d-1})$ such that $||g||_k < 4\epsilon _0$ and $||g^{N_0}-F||_k < r$. $\square $  
 
 \bigskip
 
 As for the perfectness results, we quote the following theorem from the literature. We have used only the perfectness claim in this theorem (but not the uniform perfectness with the given strong bound) which has been proved in a number of sources mentioned in the introduction.
 
 \begin{lem} \label{lem:hrt} Let $M$ be a compact smooth manifold. Then the group $\mathrm{Diff}_c^{\infty }(M)_0$ is perfect. Moreover, there exists an open neighborhood $U$ of identity such that for all $f\in \mathrm{Diff}_c^{\infty }(M)$ there exist $u_1, \dots , u_8\in \mathrm{Diff}_c^{\infty }(M)$ such that $f = [u_1, u_2][u_3,u_4][u_5,u_6][u_7,u_8]$.
 \end{lem}
 
 \medskip
 
 {\bf Proof.} See \cite{HRT}.
 
 \bigskip
 
 {\em Acknowledgement:} I am thankful to Michael Cohen for stimulating discussions.


\bigskip
 
 \begin{thebibliography}{99}
 
 \bibitem{A} A.Akhmedov, Questions ans remarks on discrete and dense subgroups of Diff(I). {\em Journal of Topology and Analysis},  vol. 6, no. 4, (2014), 557-571
 
 \bibitem{AC} A.Akhmedov and M.P.Cohen, Existence and genericity of finite topological generating sets for homeomorphism groups. Preprint. http://arxiv.org/abs/1508.04604
   
 \bibitem{Ep} D.B.A.Epstein, The simplicity of certain groups of homeomorphisms, Composito Math {\bf 22}, (1970), 165-173.
 
 \bibitem{GL} T.Gelander  and F.Le Maitre, Local topological generators and quasi non-archimedean topological groups. Preprint. http://arxiv.org/abs/1505.00415
 
  \bibitem{GS} \'E. Ghys and  R. Sergiescu, Sur un groupe remarquable de diffeomorphismes du cercle, {\em Comment. Math. Helvetici} 62 (1987) 185-239.
  
  \bibitem{Gr} R. Grzaslewicz, Density theorems for measurable transformations, {\em Colloq. Math.} {\bf (48) 2}, 245-250, 1984.
  
  \bibitem{He1} M.R. Herman, Simplicit\'e du groupe des diff\'eomorphismes de classe $C^{\infty }$, isotopes \'a l'identit\'e, du tore de dimesion $n$, {\em C. R. Acad. Sci. Paris Sr. A} {\bf 273} (1971) 232-234.
  
  \bibitem{He2} M.R. Herman, Sur le groupe des diff\'eomorphismes du tore. {\em Ann. Inst. Fourier (Grenoble)}, {\bf 23} (1973), 75-86.

   
 \bibitem{HM} K. H. Hofmann and S. A. Morris, Weight and $c$, {\em Journal of Pure and Applied Algebra}, vol. 68, no. 1--2 (1990), 181-194.
 
\bibitem{HRT} S.Haller, T.Rybicki and J.Teichmann. Smooth perfectness for the group of diffeomorphisms. {\em J.Geom. Mech.} {\bf 5}, (2013) 281-294  
 
\bibitem{K} M. Kuranishi, Two elements generations on semi-simple Lie groups.  {\em Kodai. Math. Semin. Rep.} (1949), 89-90.
 
\bibitem{Ka} M.I. Kadets, Proof of topological equivalence of all separable infinite dimensional Banach spaces, {\em Soviet Math. Dokl.}, 7 (1966), 319-322.   

\bibitem{KR} A. S. Kechris and C. Rosendal, Turbulence, amalgamation and generic automorphisms of homogeneous structures, {\em Proc. London Math. Soc.} 94 (2007), no. 2, 302-350.
 
\bibitem{Ma1} J.N.Mather, Commutators of diffeomorphisms, {\em Comment. Math. Helv.} {\bf 49} (1974), 512-528.

\bibitem{Ma2} J.N.Mather, Commutators of diffeomorphisms II, {\em Comment. Math. Helv.} {\bf 50} (1975), 33-40.

\bibitem{Mac} D.Macpherson, Groups of automorphisms of $\aleph _0$-categorical structures. {\em Quart. J. Math. Oxford} {\bf 37 (2)}, 449-465, 1986. 

\bibitem{P} V.S.Prasad, Generating dense subgroups of measure preserving transformations, {\em Proc. Amer. Math. Soc.}, {\bf (83) 2}, 286-288, 1981.
 
\bibitem{T} T. Tsuboi, On the perfectness of groups of diffeomorphisms of the interval tangent to the identity at the endpoints. {\em Proceedings of Foliations, Geometry and Dynamics}, Warsaw, 2000. ed. by Pawel Walczak et al. 

\bibitem{Th} W. Thurston, Foliations and groups of diffeomorphisms, {\em Bull. Amer. Math. Soc.} {\bf 80} (1974) 304-307.
  
 \end{thebibliography}
 \end{document}